\documentclass[12pt]{article}
\usepackage{textcomp, amsmath, amssymb, amsthm, enumerate}
\usepackage{graphicx}
\usepackage{yhmath}

\newtheorem{theorem}{Theorem}[section]
\newtheorem{lemma}[theorem]{Lemma}

\newtheorem{remark}[theorem]{Remark}

\def\rr{{\mathbb R}}

\def\sik{{\rr}^2}

\def\am{^{-1}}

\def\cc{{\mathbb C}}

\def\su{\subset}

\def\se{\setminus}

\def\al{\alpha}
\def\be{\beta}
\def\ga{\gamma}
\def\Ga{\Gamma}
\def\de{\delta}

\def\ep{\varepsilon}

\def\la{\lambda}

\def\cd{\cdot}
\def\stb{,\ldots ,}

\def\emp{\emptyset}
\def\del{\partial}

\def\nl{[0,1]}

\def\V{\Vert}
\def\msk{\medskip}
\def\bsk{\bigskip}
\def\noi{\noindent}

\def\ol{\overline}

\def\sumin{\sum_{i=1}^n}

\def\sumj0m{\sum_{j=0}^m}

\def\sumi0n{\sum_{i=0}^n}

\def\proof{\noi {\bf Proof.} }

\def\x1n{x_1 \stb x_n}
\def\y1n{y_1 \stb y_n}

\def\cl{{\rm cl}\, }

\def\dist{{\rm dist}\, }
\def\Int{{\rm Int}\, }
\def\ext{{\rm ext}\, }

\date{August 7, 2016}

\begin{document}
\title{Closed sets with the Kakeya property}
\author{M. Cs\"ornyei, K. H\'era and M. Laczkovich}

\footnotetext[2]{The research of the third author was partially
supported by the Hungarian National Foundation for Scientific Research, 
Grant No. K104178}

\maketitle

\begin{abstract}
We say that a planar set $A$ has the Kakeya property if 
there exist two different positions of $A$ such that $A$
can be continuously moved from the first position to
the second within a set of arbitrarily small area.
We prove that if $A$ is closed and has the Kakeya property,
then the union of the nontrivial connected
components of $A$ can be covered by a null set
which is either the union of parallel lines or the union of concentric
circles. In particular, if $A$ is closed, connected
and has the Kakeya property,
then $A$ can be covered by a line or a circle.
\end{abstract}

\section{Introduction and main results}\label{s1}
It is well-known
that a line segment can be continuously moved in a planar
set of arbitrarily small area such 
that it returns to its starting position,
with its direction reversed. 
This fact was first proved by Besicovitch as a solution
to the classical Kakeya problem \cite{B}.

We may ask if there are other sets in the plane having a similar
property. We shall say that a planar set $A$ has the {\it Kakeya property},
shortly property (K), if there exist two different positions of 
$A$ such that one can move $A$ continuously from the 
first position to the second
such that the set of points touched by the moving set has 
arbitrarily small area.
Obviously, every line segment has property (K). More generally,
if $A$ is a null set, and $A$ is the union of parallel lines, then
$A$ has property (K). Indeed, sliding $A$ in the direction of the lines
we only cover a set of measure zero. Similarly, 
if $A$ is a null set, and $A$ is the union of concentric circles,
then $A$ has property (K), as rotating $A$ about the center of the 
circles we only cover a set of measure zero. 
We shall say that a set $A\su \sik$ is a {\it trivial} (K)-set, if $A$
can be covered by a null set 
which is either the union of parallel lines or the union of concentric
circles. 
Our aim is to prove that 
if $A$ is closed and has the Kakeya property, then 
the union of the nontrivial connected components of $A$ is trivial.

Besicovitch's construction used the `P\'al joins' in order to shift
the line segment to an arbitrary parallel position using 
arbitrarily small area. (There are solutions to the Kakeya
problem without using P\'al joins, see e.g. \cite{C1}, \cite{C2}, \cite{D}).
Using the P\'al joins and Besicovitch's theorem 
it is clear that every line segment
can be moved to arrive at any prescribed position
within a set of arbitrarily small area.

We shall say that 
a planar set $A$ has the {\it strong Kakeya property}, shortly
property ${\rm (K^s )}$ if, for every prescribed position, $A$ 
can be continuously moved within a set of 
arbitrarily small area to the prescribed final position.

The question, whether or not circular arcs have property ${\rm (K^s )}$
was asked by F. Cunningham (see \cite[p.~591]{C2}). 
It is proved in
\cite{HL} and \cite{CCL} that the answer is affirmative, at least for
circular arcs of angle short enough.
In this note we shall prove that
every connected closed set with property ${\rm (K^s )}$ 
is a line segment or a circular arc.

Now we formulate these notions and results precisely. A {\it rigid motion} is
an isometry of the plane preserving orientation; that is a translation
or a rotation. We denote the identity map by $j$. 
We identify the plane with the complex plane $\cc$, and put
$S^1 =\{ x\in \cc \colon |x|=1\}$.
Then every rigid motion
is a function of the form $x\mapsto ux+c$, where $u,c\in \cc$ and 
$u\in S^1$.
By a {\it continuous movement} $M$ we mean a 
map $t \mapsto M_t$ $(t\in [0,1])$ such that $M_t$ is a rigid motion 
for every $t \in [0,1]$, the map $(t,x) \mapsto M_t (x)$
is continuous on $\nl  \times \mathbb{R}^2$, and $M_0 =j$.
If $M$ is a continuous movement, then the set of points 
touched by a moving set $A$ is
$$W_M (A)=\{ M_t (x) \colon t\in \nl , \ x\in A\} .$$
The two dimensional Lebesgue measure of the set $A\su \sik$ is denoted
by $\la (A)$. 

The set $A$ has property (K) if there is a rigid motion $\al \ne j$
such that for every $\ep >0$
there exists a continuous movement $M$ such that $M_1 =\al$ and $\la (W_M (A))
<\ep$. The set $A$ has property ${\rm (K^s )}$ if for every rigid motion $\al$
and for every $\ep >0$
there exists a continuous movement $M$ such that $M_1 =\al$ and $\la (W_M (A))
<\ep$. By a nontrivial connected component of a set $A$ we mean
a connected component of $A$ having at least two points. Our main 
results are the following.

\begin{theorem} \label{t1}
Let $A\su \sik$ be a closed set having
property {\rm (K)}. Then the union of the nontrivial connected components of 
$A$ is a trivial {\rm (K)}-set.  If $A$ is nonempty, closed, connected and has 
property {\rm (K)}, then $A$ is a line segment, a half line, a line,
a circular arc, a circle or a singleton.
\end{theorem}

\begin{theorem} \label{t2}
If $A\su \sik$ is a nonempty closed and connected set having
property {\rm ${\rm (K^s )}$}, then
$A$ is a line segment, a circular arc or a singleton.
\end{theorem}

It is easy to see that full circles  do not have property ${\rm (K^s )}$.
Indeed, moving a circle $C$ to a circle disjoint from $C$, the moving
circle must touch every point inside $C$.
Since every circle has property (K),
we can see that properties (K) and ${\rm (K^s )}$ are not equivalent. 
As we mentioned earlier, all circular arcs of angle short enough 
have property ${\rm (K^s )}$.
It is not clear, however, whether or not all circular arcs other than
the full circles have property ${\rm (K^s )}$, while,
obviously, they have property (K).

It is easy to see that {\it every set of 
Hausdorff dimension less than $1$ has property} ${\rm (K^s )}$.
We show that {\it there are compact sets with property ${\rm (K^s )}$ having  
Hausdorff dimension $2$.} Let $E\su [1,2]$ and $F\su [0,1/2]$ be compact 
sets of linear measure zero
and of Hausdorff dimension $1$. We put 
$$A=\{ (x,y)\in \sik \colon x>0, \ y\in F,\ \sqrt{x^2 +y^2} \in E\} .$$
It is easy to check that $A$ is a bi-Lipschitz image of $E\times F$,
and thus $A$ has Hausdorff dimension $2$. We prove that $A$ has property ${\rm (K^s )}$. 
Let $\al (x)=ux+c$ be a rigid motion, where $|u|=1$. Let $c=r\cd v$,
where $r\ge 0$ and $|v|=1$. Then the following
continuous movement brings $A$ onto $\al (A)$. First we rotate $A$
about the origin to obtain $v\cd A$, then we translate $v\cd A$ by
the vector $c$ to obtain $v\cd A+c$, and then we rotate $v\cd A+c$ about
the point $c$ to obtain $u\cd A+c =\al (A)$. If
$B=\{ (x,y)\colon \sqrt{x^2 +y^2} \in E\}$, then the set of points 
touched by this
continuous movement is a subset of 
$$B\cup v\cd \{ (x,y)\colon y\in F\} \cup (B+c).$$
This is a set of plane measure zero, showing that $A$ 
has property ${\rm (K^s )}$.

We note that the union of finitely many parallel 
segments have property ${\rm (K^s )}$;
this was proved by Roy O. Davies in \cite[Theorem 3]{D}. It is not
clear if the union of finitely many
concentric circular arcs can have property ${\rm (K^s )}$.

\begin{remark}
There are compact subsets of the plane having property {\rm (K)}
without being trivial {\rm (K)}-sets.
\end{remark}
We sketch the construction of a totally disconnected compact 
set $A$ having property {\rm (K)} 
such that $A$ has orthogonal projection of positive linear measure 
in every direction $e \in S^1$, and the distance set
$$D(A;p)=\{|x-p| \in \rr \colon x \in A \}$$ 
has positive 
linear measure for all $p \in \sik$. This easily implies that $A$ is not 
trivial. The construction is based on the so-called 
'venetian blind' construction.
Using this method M. Talagrand proved in \cite{T}  that for 
every u.s.c. function $f$ on
$S^1$ there exists a compact set $K \su \sik$ such that the measure 
of the projection of $K$
in direction $e$ is $f(e)$ for every $e \in S^1$. 
We say that a rectangle $R$ is in the direction $e$ if the longer side
of $R$ is of direction $e$.

We fix a direction $e$ and choose a convergent sequence 
of distinct directions $e_n \to e$. 
We shall construct a decreasing sequence of compact sets $K_n$ such that
$K_n$ is the union of finitely many narrow rectangles in the direction $e_n$ for
every $n$.
Let $K_1$ be a closed rectangle in the direction $e_1$ of width less than $1$. 
Then the measure of the projection of $K_1$ in the direction $e_1$ is 
less than $1$, and 
the measure of the distance set $D(K_1 ;p)$
is positive for all $p\in \sik$.
In the $(n+1)^{\rm th}$ step we replace each member of the finite system of 
closed rectangles whose union is
$K_n$ by a finite union of closed rectangles such that the new rectangles
have direction $e_{n+1}$, and 
the union of these $(n+1)^{\rm th}$ generation rectangles,
$K_{n+1}$, is a subset of $K_n$. Moreover, by choosing these rectangles 
appropriately, 
the measure of the projection of $K_{n+1}$ in the direction $e_{n+1}$ will be
less than $1/(n+1)$, and the measure of the projection of $K_{n+1}$ in 
directions $e_1,\dots,e_n$ will be only slightly smaller than that of
the projection of $K_n$ in directions $e_1 ,\dots , e_n$, respectively.
If we choose the new rectangles originating from a fixed 
previous rectangle to be close enough 
to each other, we can achieve that the set $D(K_{n+1} ;p)$
is only slightly smaller than $D(K_{n} ;p)$
for all $p \in \sik$, $|p| \le n$. We set $K=\bigcap_{n=1}^{\infty} K_n$. 
We may also ensure that the projection of $K$ 
has positive linear measure in all directions (see \cite{T}). 
The set $D(K ;p)$ has positive measure for all 
$p \in \sik$. Indeed, for every $p$ 
there is an $n$ such that $|p| \le n$ and then,
for all $m \ge n$ we have 
$\la_1( D(K_{m} ;p)) \ge \la_1( D(K_{n} ;p)) /2 $
and thus $\la_1( D(K ;p)) \ge  \la_1( D(K_{n} ;p)) /2 >0$.
It is easy to see that $K$ has 
property {\rm (K)}. Indeed, if $\al$ is a translation 
by a vector of direction $e$ and length $R>0$, then $K$ can be moved continuously 
to $\al(K)$ within a set of arbitrarily small area. 
Then $K$ has the desired properties, and the proof is finished.

\section{Some auxiliary results and the proof of Theorems 
\ref{t1} and \ref{t2}} \label{s2}

In this section we formulate some auxiliary results needed
for the proof of Theorems \ref{t1} and \ref{t2}.

The translation by the vector $c$ will be denoted by $T_c$. Thus
$T_c (x)=x+c$ for every $x\in \cc$. The rotation about the point
$c$ by angle $\phi$ is denoted by $R_{c, \phi}$. Thus
$R_{c, \phi} (x)=e^{i\phi} (x-c)+c$ for every $x\in \cc$. 
If $\al$ is a rigid motion and $\al ^2 \ne j$, then we define the 
\emph{elementary movement} $E^\al$ determined by $\al$ 
as follows. If $\al =T_c$, then we put
$E^\al_t =T_{tc}$ for every $t\in \nl$. If $\al =R_{a, \phi}$, then
the condition $\al ^2 \ne j$ implies that $\phi \not\equiv \pi$
(mod $2\pi$). Therefore, we may assume that $|\phi |<\pi$.
Then we define $E^\al_t =R_{a, t\phi}$ for every $t\in \nl$.

The inverse of the rigid motion $\al$ is denoted by $\al \am$. It is clear
that if $M$ is a continuous movement, then so is $M\am$, where
we put $M\am _t =(M_t )\am $ $(t\in \nl )$.
The $\ep$-neighbourhood of a set $E\su \sik$ is defined by
$$U(E,\ep )=\{ x\in \sik \colon \dist(x,E)<\ep \} .$$

We denote by $L$ the space of functions
$x\mapsto ux+v$ $(x\in \cc )$, where $u,v\in \cc$.
Clearly, $L$ is a linear space over $\cc$ with pointwise addition and 
multiplication by constants. We endow $L$ with the norm
$\V f \V=|u| + |v|$. It is easy to check that 
$\V f \V  =\sup \{ |f(x)| \colon x\in \cc, \ |x|\le 1\}$
for every $f\in L$, and thus $\V \cd 
\V$ is indeed a norm.

It is also easy to check that if $f_1(x)= u_1 x+v_1$ and 
$f_2(x)=u_2 x+v_2$ are rigid motions, then 
\begin{equation}\label{e8}
\V f_1^{-1} - f_2^{-1} \V \le (1+|v_2|) \V f_1-f_2 \V.
\end{equation}
\begin{lemma}\label{l3}
If $A\su \sik$ has property {\rm (K)}, then there
exists a rigid motion $\al$ such that $\al ^2 \ne j$,
and the following condition is satisfied.
For every $\ep >0$ there is a
continuous movement $M$ such that
$M_1 =\al$, $\la (W_{M} (A))<\ep$, and $ \V M_t - E_t ^\al  \V < \ep$ 
for every $t \in [0,1]$.
\end{lemma}
By a continuum we mean a compact connected set.
The set $A\su \sik$ is said to
be irreducible between the distinct points $a$ and $b$ provided that 
$A$ is connected, $a,b\in A$,
and these two points cannot be joined by any closed, connected, proper
subset of $A$. We shall need the following topological lemma on irreducible
continuums.
\begin{lemma}\label{l4}
Let $A\su \sik$ be a continuum which is 
irreducible between the distinct points $a$ and $b$, and suppose that
$\sik \se A$ is connected. Let $D$ be an
open disc not containing the points $a$ and $b$. 
Then every neighbourhood of every point of $A\cap D$ intersects at least
two of the connected components of $D\se A$.
\end{lemma}
The proof of Lemmas \ref{l3} and \ref{l4} will be given in the next
two sections. 
\begin{lemma}\label{l5}
Let $A\su D\su \sik$ be arbitrary and $G\su D\se A$.
Suppose that $M$ is a continuous movement, $t\in \nl$, and
$M\am _s (x)\in D$ for every $s\in [0,t]$ and $x\in G$.
If $G$ and $M\am _t (G)$ are subsets of
distinct connected components of $D\se A$, then $G\su W_M (A)$.
\end{lemma}
\proof 
Let $u\in G$ be arbitrary. Clearly, the map $\ga \colon [0,t] \to D$
defined by $\ga (s) = M\am  _s (u)$ $(s\in [0,t])$ is a continuous curve. Now 
$\ga (0)=u\in G$ and $\ga (t)= M\am _t (u)\in M\am _t (G)$.
Since $G$ and $M\am _t (G)$ are subsets of
distinct connected components of $D\se A$ and $\ga ([0,t])\su D$, 
it follows that there exists an $s \in [0,t]$ such that 
$\ga (s)\in A$. If $\ga (s)=a$, then $M\am  _s (u) =a$ and
$u=M_s (a)\in M_s (A)\su W_M (A)$. This is true for every
$u\in G$, and thus $G\su W_M (A)$.
\hfill $\square$

\msk \noi
{\bf Proof of Theorem \ref{t1}} (subject to Lemmas \ref{l3} and \ref{l4}).
We denote by $B(x,r)$ the open disc with centre $x$ and radius $r$.
Let $A\su \sik$ be a closed set having
property {\rm (K)}. By Lemma \ref{l3}, there exists 
a rigid motion $\al$ satisfying the conditions of
Lemma \ref{l3}.  Let $A'$ denote the union
of all nontrivial components of $A$. Since $A'\su A$
and $A$ has property (K), we have $\la (A)=0$ and $\la (A')=0$.
We shall prove that if $\al$ is a 
translation by the vector $v\ne 0$, then every nontrivial connected 
component of $A$ is covered by a line parallel to $v$, and if 
$\al$ is a rotation about the point $c$, then every nontrivial connected 
component of $A$ is covered by a circle of centre $c$. Clearly, since $A'$ has measure zero, therefore it cannot meet positively many parallel lines in positive length, and similarly, it cannot meet positively many concentric circles in positive length. Therefore indeed $A'$ is a trivial (K)-set.

\medskip

%Then we complete the proof by the following argument. 
%Suppose that the nontrivial components of $A$
%can be covered by parallel lines. Translate $A'$ by every rational 
%distance in the direction of these lines, and let $A''$ be the union of these
%translated copies. Then $\la (A'')=0$, and $A''$ is the union of parallel lines 
%which cover $A'$. We argue similarly if the nontrivial components of $A$
%can be covered by concentric circles. 
%It is clear that $A'$ is a trivial (K)-set.

We shall only prove the statement when $\al$ is a rotation; the
case when $\al$ is a translation can be treated similarly.

Let $\al$ be a rotation. We may assume that the centre of rotation is
the origin. Let $A_1$ be a connected component of $A$. We
have to prove that $A_1$ is covered by a circle of centre $0$.
Suppose this is not true. 
Then the set $\Ga =\{ |x| \colon x\in A_1 \}$ is a nondegenerate interval. 
Let $r_1 ,r_2 \in \Ga$ be such that $0<r_1 <r_2$.
We may assume that $r_2 <1$, since otherwise we replace $A$ by a suitable
similar copy.

Let $U$ denote the annulus $\{ x \colon r_1 <|x|<r_2 \} \su B(0,1)$. 
The set $A_1$ contains an irreducible (that is, minimal) connected closed subset
$C$ such that it intersects both of the circles $|x|=r_1$ and $|x|=r_2$.
(The existence of such a set follows from 
\cite[Theorem 2, \S 42, IV, p. 54]{K}.) 
Then $C$ is contained in $\cl U$. Indeed, every connected component of
$C\cap \cl U$ is a quasi-component of $C\cap \cl U$ by 
\cite[Theorem 2, \S 47, II, p. 169]{K}. Using this fact one can
prove that there is a connected component $C_1$ of $C\cap \cl U$ which intersects 
the circles $|x|=r_1$ and $|x|=r_2$. Thus, by the minimality of $C$, we have
$C_1 =C$ and $C\su \cl U$.

Let $C^*$ denote the set of points $p\in C\cap U$
with the following property: if $B(p,r)\su U$, then 
every neighbourhood of $p$ intersects at least two connected components of $B(p,r)\se C$.
First we prove that if $p\in C^*$, then $C$ contains a subarc of the circle
$\{ x\colon |x|=|p|\}$. 

Let $p\in C^*$, and fix an $r>0$ with $B(p,r)\su U$. Put $D=B(p,r)$.
There is a $0<t_0 \le 1$
such that the arc $I=\{ E^{\al \am}_t (p)\colon t\in [0,t_0 ]\}$ is in $D$.
We prove that $I\su C$. Suppose this is not true. Then there is a 
$0<t_1 \le t_0$ such that $q=E^{\al \am}_{t_1} (p)\in D\se C$. Let 
$\de >0$ be such that
$B(q,\de )\su D\se C$. Since $p\in C^*$, the neighbourhood
$B(p,\de )$ intersects at least two connected components of $D\se C$. Since 
$B(q,\de )$ belongs to a connected component of $D\se C$, there is an
open disc $G$ such that $\cl G\su B(p,\de )\se C$, and 
$G$ and $B(q,\de )$ belong to distinct connected components of $D\se C$.

Since $\al$ satisfies the conditions of 
Lemma \ref{l3}, it follows that for every $\ep >0$ there is a
continuous movement $M$ such that
$\la (W_{M} (A))<\ep$ and $\V M_s  -E^{ \al}_s \V < \ep $
for every $s \in [0,t_1]$. Therefore, applying \eqref{e8} with
$f_1 =M_s$ and $f_2 =\al$ we obtain
$\V  M_s \am -E_s ^{\al \am} \V <\ep$ 
for every $s \in [0,t_1]$. (Note that we have $v_2 =0$ in this case.) 
Since $G \su U \su B(0,1)$, we obtain
$$|M_s\am (x) -E^{ \al \am}_s (x)| \le \V  M_s\am -E_s ^{\al \am} \V  <\ep$$
for all $x \in G$ and $s \in [0,t_1]$.
It is clear that if $\ep$ is small enough,
then $M\am _s (G)\su U(E^{ \al \am}_s (G), \ep) \su D$ for every $s\in [0,t_1]$. 
Also, we have $E^{\al \am} _{t_1} (p)=q$. Since
$E^{\al \am} _{t_1}$ is a rigid motion, $\cl E^{\al \am} _{t_1} (G)\su B(q,\de )$.
Therefore, for small $\ep$,
we have $M\am _{t_1} (G)\su B(q,\de )$. Thus, by Lemma \ref{l5},
we have $G\su W_{M} (C)$. Thus $\la (G)\le \la (W_M (C))$.
Since $\la (W_M (C))<\ep$, this is impossible if $\ep <\la (G)$.
This contradiction shows that $I\su C$ as we stated.

There are points $a,b\in C$ such that $|a|=r_1$ and $|b|=r_2$.
It is clear that $C$ is irreducible 
between $a$ and $b$. Suppose that $\sik \se C$ is connected.
Then, by Lemma \ref{l4}, we have $C^* =C\cap U$. Then, by what we proved above,
every point of $C\cap U$ is covered by a circular arc of centre $0$
belonging to $C$.
Since $C$ is connected, this implies that for every $r_1 <r<r_2$, the
circle $|x|=r$ contains an arc belonging to $C$. This, however, is impossible,
since $\la (C)=0$.

Therefore, the set $\sik \se C$ cannot be connected. Let $V$ be a
bounded connected component of $\sik \se C$. Clearly, we have $V\su B(0, r_2 )$.  
We show that $C$ contains a full circle of centre $0$.
This is clear if $V= B(0, r_1 )$ since, in that case, $\del V$ is a
circle and is contained in $C$. So we may assume that $V\ne B(0,r_1 )$.
Then $V\cap U\ne \emp$.
We prove that $\cl V$ is an annulus of centre $0$.
Suppose this is not true. Then there are points $x_1 \in \cl V$
and $x_2 \notin \cl V$ such that $|x_1 |=|x_2 |$.
Then $B(x_2 ,\ep )\cap V =\emp$ for $\ep$ small enough. Since
$B(x_1 ,\ep )\cap V \ne \emp$ and $C$ is nowhere dense, it is clear that
there are points $y_1 \in V$
and $y_2 \in B(x_2 ,\ep )$ such that $|y_1 |=|y_2 |$.
Then $y_2$ belongs to a connected component of $\sik \se C$ different from $V$.

This easily implies that there is an $\eta >0$ such
that for every $|y_1 | -\eta <r<|y_1 |$ there is a point $p\in \del (\cl V)$ 
with $|p|=r$. It is easy to see that if $p\in U\cap \del (\cl V)$, then
$p\in C^*$. Indeed, let $D$ be an open disc such that $p\in D\su U$.
Then $D$ intersects $V$, and thus it intersects 
at least two connected components of $U\se C$, since otherwise $D$
would be a subset of $\cl V$, contradicting $p\in \del (\cl V)$.
This is true for every disc $B(p,\de )\su D$ as well, and thus 
$B(p,\de )$ intersects at least two connected components of $D\se C$, proving
that $p\in C^*$. As we saw above, this implies that for every 
$|y_1 | -\eta <r<|y_1 |$, $C$ contains a subarc of the circle $|x|=r$.
This, however, contradicts the fact that $C$ is of measure zero. 

Therefore, the set $\cl V$ must be an annulus. Since 
$\del (\cl V)\su \del V\su C$, it follows that $C\cap \cl U$ 
contains a full circle of centre $0$. 

We proved that $A_1 \cap \{ x\colon r_1 <|x|<r_2 \}$ 
contains a full circle of centre $0$ for every $r_1 ,r_2 \in \Ga$
with $r_1 <r_2 <1$. Thus $A_1$ contains a dense subset of 
$\{ x\colon |x|\in \Ga \} \cap B(0,1)$.
Since $A_1$ is closed, we find that $A_1$ must contain the whole
set $\{ x\colon |x|\in \Ga \} \cap B(0,1)$, which is clearly impossible.
Thus Theorem \ref{t1} is proved.

In order to prove Theorem \ref{t2} we only have to show that
halflines and full circles do not have property 
${\rm (K^s )}$. We already saw this for circles.
The case of halflines is also clear, as a horizontal
halfline cannot be moved continuously to a vertical halfline
touching only a finite area. 

\section{Proof of Lemma \ref{l3}}\label{s3}
Let $\al$ be a rigid motion and $\al (x)=u x+a$ where $|u|=1$. 
Let $\al ^n$ denote the $n^{\rm th}$ iterate of $\al$.
We would like to compare the distances $\V \al^n-j \V$ and $\V \al-j \V$. 
 
It is easy to check that if $\al (x)=u x+a$ where $|u|=1$, then
$\al ^n (x)=u^n x+(1+u+\ldots +u^{n-1})a$, and thus
\begin{equation}\label{e1}
\begin{split}
\V \al ^n-j \V &=|u^n -1|+|1+u+\ldots +u^{n-1}|\cd |a|=\\
&=|1+u+\ldots +u^{n-1}|\cd (|u-1|+|a|)=\\
&=|1+u+\ldots +u^{n-1}|\cd \V \al-j \V .
\end{split}
\end{equation}
Since $|u|=1$, it follows from \eqref{e1} that 
$\V \al ^n -j \V \le n\cd \V \al -j\V$ for every $\al$. 
Now we show that if $|u-1|\le 1/n$, then
\begin{equation}\label{e2}
\V \al ^n -j\V \ge \frac{n}{2} \cd \V \al -j\V .
\end{equation}
Since $|u^i -u^{i-1}|=|u-1|$ for every $i=1,2,\ldots$ we obtain
$|u^k -1|\le k\cd |u-1|$ for every $k=1,2,\ldots$. Then,
assuming $|u-1|\le 1/n$ we find that
$$|(1+u+\ldots +u^{n-1}) -n| \le (1+\ldots +(n-1))\cd |u-1| \le 
\frac{n(n-1)}{2}\cd  \frac{1}{n} <\frac{n}{2} .$$
Thus $|1+u+\ldots +u^{n-1}| \ge n/2$, and then \eqref{e1} gives \eqref{e2}.

It is clear that if $\al$ is a translation or a rotation of angle $\phi$
with $|\phi |<\pi /n$, then
\begin{equation}\label{e3}
E^{\al ^n}_{i/n} =\al ^i \qquad (i=1\stb n).
\end{equation}
We shall also need the following estimate.
\begin{lemma}\label{l1}
Let $\al$ be a rigid motion. If $\al ^2 \ne j$ then, for every 
$t_1 ,t_2 \in \nl$ and $|x|\le 1$ we have
\begin{equation}\label{e5}
\left| E^\al _{t_1} (x)-E^\al _{t_2} (x) \right| \le 2|t_1 -t_2 |\cd 
(\V \al \V +1).
\end{equation}
\end{lemma}

\proof
If $\al =T_c$, then 
$$E^\al _{t_1} (x)-E^\al _{t_2} (x)= (x+t_1 c)-(x+t_2 c) =(t_1 -t_2 )c,$$
from which \eqref{e5} is clear. 

Now let $\al =R_{c,\phi}, \
R_{c,\phi}(x) = e^{i\phi} (x-c)+c= e^{i\phi} x +c (1-e^{i\phi})$, where $|\phi | < \pi$. Then
$$\V \al \V =|e^{i\phi}| +|c|\cd |1-e^{i\phi}|=
1+|c|\cd |1-e^{i\phi}|= 1 + 2 \cd|\sin (\phi /2)|\cd |c|.$$
Let $S$ denote the circle having centre $c$ and radius $|x-c|$. 
Then $E^\al _{t_1} (x)$ and $E^\al _{t_2} (x)$ are the endpoints of a
subarc of $S$ of angle $t_1 \phi -t_2 \phi$. Therefore, assuming $|x|\le 1$
and putting $h=|t_1 -t_2 |$, we have
\begin{equation}\label{e6}
\left| E^\al _{t_1} (x)-E^\al _{t_2} (x)\right| \le
|t_1 -t_2 |\cd |\phi | \cd |x-c|=h\cd |\phi | \cd |x-c| \le
h\cd |\phi | \cd (1 +|c|).
\end{equation}
Since $\sin x\ge (2/\pi )\cd x$ for every $x\in [0,\pi /2]$,
we have $|\sin (\phi /2)|\ge |\phi /\pi |$ by $|\phi /2|< \pi /2$.
Thus, by \eqref{e6} we get
\begin{align*}
\left| E^\al _{t_1} (x)-E^\al _{t_2} (x)\right| & \le 
   h\cd  \pi \cd |\sin (\phi/2)| \cd (1+|c|) = \\
 &= h\cd  (\pi/2) \cd 2|\sin (\phi/2)| \cd (1+|c|) = \\
 &= h\cd  (\pi/2) \cd 2|\sin(\phi/2)| \cd |c| +h \cd (\pi/2) \cd 2|\sin (\phi/2)| \le \\
& \le h\cd  (\pi/2) \cd 2|\sin (\phi/2)| \cd |c| + h \cd (\pi/2) \cd 2 = \\
&= h\cd  (\pi/2) ( \V \al \V +1) \le 2h ( \V \al \V +1),
\end{align*}
proving \eqref{e5}. \hfill $\square$

\begin{lemma}\label{l2}
Let $\al _n (x)=u_n x+v_n$ $(n=1,2,\ldots )$ and $\al (x)=ux+v$
be rigid motions, 
where $u_n \to u$ and $v_n \to v$ as $n\to \infty$. 
Suppose that $\al ^2 \ne j$.
Then we have $E^{\al _n}_t (x)\to E^{\al}_t (x)$ 
uniformly on the set
$\{ (t,x)\colon t\in \nl , \ |x|\le 1\}$.
\end{lemma}
In order to prove the lemma one has to distinguish between the
cases $u=1$ and $u\ne 1$. In both cases the proof is an easy 
computation; we leave the details to the reader.

Now we turn to the proof of Lemma \ref{l3}.
Since $A\su \sik$ has property (K), there is a rigid motion $\be \ne j$
such that for every $\ep >0$
there exists a continuous movement $M$ with $M_1 =\be$ and $\la (W_M (A))
<\ep$. Let $\be (x)=v_0 x+b_0$, where $|v_0 |=1$.
We put $\eta =\min (|v_0 -1|+|b_0 | ,1)$. Since $\be \ne j$,
we have $\eta >0$.

We choose, for every $n$, a 
continuous movement $M^n$ such that $M^n_1 =\be$ and $\la (W_{M^n} (A))
<1/n^2$. Let $M^n_t  (x)=v_n (t)x+b_n(t)$ for every $x\in \cc$
and $t\in \nl$, where $|v_n |=1$.
Since the maps $t\mapsto v_n (t)$, $t\mapsto b_n(t)$ are continuous
and $v_n (1)=v_0 , \ b_n(1)=b_0$, there is a $0<t_n \le 1$ such that
$|v_n (t_n )-1|+|b_n(t_n )|=\eta /n$, and 
$|v_n (t )-1|+|b_n(t )|<\eta  /n$ for every $t\in [0,t_n )$.
Replacing $M^n$ by $M^n \circ \psi$, where $\psi$ is a 
suitable homeomorphism of $\nl$ onto itself, we may assume that $t_n =1/n$.

We put $\be _n =M^n_{1/n}$ and $\al _n =\be _n^n$.
Then $\V \be _n -j\V =\eta /n$ and
$\eta /2 \le \V \al _n -j\V \le \eta$ by \eqref{e2}.
Let $\al_n (x) =u_n x+a_n$, where $|u_n |=1$. Since
$\eta /2 \le |u_n -1|+|a_n |\le \eta$ for every $n$,
the sequences $(u_n )$ and $(a_n )$ have convergent subsequences.
Turning to a suitable subsequence we may
assume that the sequences $(u_n )$ and $(a_n )$ converge. Let
$u_n \to u$ and $a_n \to a$. Then $|u|=1$, and thus $\al (x)=ux+a$
defines a rigid motion. Our aim is to show that $\al$ satisfies the
requirements. Since $\V \al-j \V =|u-1|+|a|
\ge \eta /2 >0$, we can see that $\al \ne j$. Also, by 
$\V \al -j\V \le \eta \le 1$, we have $\al ^2 \ne j$. Indeed, 
$\al ^2 =j$ would imply $u=-1$ and $\V \al-j \V \ge 2$.

We define the continuous movement $F^n$ as follows: we replace
the movement $E^{\al _n}_t$ between the moments $t=(i-1)/n$ and $t=i/n$ by 
suitable copies of the movement $M_t$ $(t\in [0,1/n])$.
More precisely, we put
$$F^n_t =\be _n^{i-1} \circ M^n _{t-(i-1)/n}$$
for $t\in [(i-1)/n,i/n]$  and $i=1\stb n$.
Since $M^n_0 =j$ and $M^n_{1/n} =\be _n$, it follows that 
$\be _n^{i-1} \circ M^n_{1/n}  = \be _n^{i} \circ M^n_0$, and 
thus the definition makes sense.
This also proves that $F^n$ is a continuous movement. Let
$W=\{ M^n_t (x)\colon x\in A, \ t\in [0,1/n]\}$. Since $W \su W_M (A)$,
we have $\la (W)<1/n^2$. Also, $W_{F^n} (A)$ is the union of
$n$ congruent copies of $W$, therefore, we have $\la(W_{F^n} (A))<1/n$.

Next we show that for every $x\in \cc$, $|x|\le 1$ and $t\in \nl$ we have
\begin{equation}\label{e4}
\left| F^n_t (x)-E^{\al _n}_t (x)\right| \le 8/n.
\end{equation}
If $t\in [0,1/n]$, then
$|v_n (t )-1|+|b_n (t )|\le \eta  /n \le 1/n$, and thus
$$|M^n_t (x)-x|=|(v_n (t) -1)x+b_n (t)|\le 
|v_n (t) -1|\cd |x|+|b_n (t)|\le (|x|+1)/n$$
for every $x\in \cc$. Let $1\le i\le n$ and 
$t\in [(i-1)/n, i/n]$ be given. Then we have 
\begin{equation}\label{e7}
\begin{split}
\left| F^n_t (x) -\be _n^{i-1} (x)\right| &= 
\left| \be _n^{i-1} \circ M^n _{t-(i-1)/n} (x) -\be _n^{i-1} (x)\right| = \\
&= \left| M^n _{t-(i-1)/n} (x) - x \right| \le (|x|+1)/n \le 2/n,
\end{split}
\end{equation} 
where we used the fact that $\be _n^{i-1}$ is an isometry. 
On the other hand, by Lemma \ref{l1} and by \eqref{e3},
we have
\begin{align*}
\left| E^{\al _n} _t (x)-\be _n^{i-1} (x)\right| &=
\left| E^{\al _n} _t (x)-E^{\al _n} _{(i-1)/n} (x) \right| \le\\
&\le (2/n) \cd (\V \al _n \V +1)\le (2/n) \cd 3 = 6/n.
\end{align*}
Note that $\V \be _n -j\V =\eta /n <1/n$ implies
that either $\be _n$ is a translation or it is a rotation by an
angle $\phi$ with $\phi |<\pi /n$, and thus \eqref{e3} can be applied.
Then, comparing with \eqref{e7}, we have \eqref{e4}.
Since $u_n \to u$ and $a_n \to a$, it follows from Lemma \ref{l2}
that $E^{\al _n}_t (x)\to E^{\al}_t (x)$ uniformly on the set
$\{ (t,x)\colon t\in \nl , \ |x|\le 1\}$. Let $\ep >0$ be given.
Then we can choose an $n$ such that $8/n<\ep /2$ and
$$\left| E^{\al _n}_t (x)- E^{\al}_t (x) \right| <\ep /2$$
for every $t\in \nl$ and $|x|\le 1$. Then, by \eqref{e4}, we find
that $|F^n_t (x)-E^{\al}_t (x)| < \ep$ for every $t\in \nl$ and $|x|\le 1$.
Thus $\V F^n_t-E^{\al}_t \V  < \ep$ for every $t\in \nl$. Since
$\la (W_{F^n} (A))<1/n<\ep$, this completes the proof.
\hfill $\square$

\section{Proof of Lemma \ref{l4}}

If $\ga \colon [a,b]\to \cc$ is a continuous curve and $p\in \cc \se
\ga ([a,b])$, then we denote by $w(\ga ;p)$ the increment of the 
argument of $\ga (t)-p$, as $t$ goes from $a$ to $b$. More precisely,
if $\phi \colon [a,b]\to \rr$ is a continuous function such that
$\phi (t)$ is one of the arguments of $\ga (t)-p$ for every $t\in [a,b]$,
then we put $w(\ga ;p)=\phi (b)-\phi (a)$. It is easy to see that
$w(\ga ;p)$ does not depend on the choice of the continuous function
$\phi$. If $\ga$ is a closed curve then $w(\ga ;p)/(2\pi )$ is 
the winding number of $\ga$ with respect to the point $p$.

Let $p$ be an arbitrary element of $A\cap D$, and let $U\su D$ 
be an open disc of centre $p$. We prove that $\cl U$ intersects at least 
two connected components of $D\se A$.

Since $A$ is irreducible between the points
$a,b\in A\se U$, it follows that $A\se U$ is not connected;
moreover, $a$ and $b$ belong to different connected components of $A\se U$.
Let $C_b$ denote the connected component of $A\se U$ containing the point $b$.

By \cite[Theorem 2, \S 47, II, p. 169]{K}, $C_b$ is a quasi-component
of $A\se U$; that is, $C_b$ is the intersection of all (relative) clopen
subsets of $A\se  U$ containing $b$. Since $a\notin C_b$, 
there is a relative clopen
set $G\su A\se U$ such that $C_b \su G$ and $a\notin G$. Then 
$(A\se U)\se G$ and $G$ are disjoint compact sets.

Since $A$ is irreducible between $a$ and $b$, it follows that $A$
is nowhere dense. Indeed, it is clear that 
otherwise it would contain proper subcontinuums containing
$a$ and $b$.

By assumption, the set $\sik \se A$ is connected. 
This implies that $\sik \se C$ is connected
for every $C\su A$. Indeed, $\sik \se C$ contains 
the connected set $\sik \se A$, and is
contained in its closure, so must be connected itself. Thus
$(A\se U)\se G$ and $G$ are disjoint compact sets such that they do not cut the
plane. Indeed, $\sik \se (A\se U)$ and $\sik \se G$
are semicontinuums, being connected open sets.

Therefore, by \cite[Theorem 9, \S 61, II, p. 514]{K},
$(A\se U)\se G$ and $G$ can be separated by a simple 
closed curve. Let $J$ be a simple
closed curve separating them. Then $J\cap (A\se U) =\emp$ and 
$J$ separates $a$ and $b$. 
By symmetry we may assume that $a\in {\rm Int} \, J$ and $b \in
{\rm Ext} \, J$. 

If $J\su D$, then $\Int J \su D$ and $a\notin \Int J$
which is impossible. Therefore, $J$ is not covered by $D$, and we 
can pick a point $a_0 \in J\se D$.

Since $A$ is connected and contains the points $a,b$, $A$ must intersect $J$.
Now we have $J\cap (A\se U)=\emp$, and thus $J\cap U\ne \emp$.
Since $a_0 \notin U$, it follows that $J\cap \del U \ne \emp$.

Let $J$ be the range of the continuous map $\ga :\nl \to \sik$, where 
$\ga (0)=\ga (1)\in \del U$, and $\ga$ is injective on $[0,1)$.
The set $F=\ga \am (\del U)$ is a closed subset of $\nl$ such that
$0,1\in F$. If $(x,y)$ is a connected component of $\nl \se F$, then $\ga ((x,y))\cap
\del U=\emp$, and thus $\ga ((x,y))$ is covered by one of $ U$ and $\ext
U$. 
By the uniform continuity of $\ga$ there exists a $\de >0$ such that
whenever $x,y\in F$ and $|x-y|<\de$, then $\ga ((x,y)) \su D$. This 
easily implies that there is a partition $0=t_0 <t_1 <\ldots <t_n =1$
such that $t_0 \stb t_n  \in F$, and for every $i=1\stb n$
either $\ga ((t_{i-1} ,t_i )) \su D$ or $\ga ((t_{i-1} ,t_i )) \cap \cl U =\emp$.

Let $J_i =\ga ([t_{i-1} ,t_i ])$ for every $i=1\stb n$. If $L_i$ denotes
the subarc of the circle $\del U$ with endpoints $\ga (t_i )$ and 
$\ga (t_{i-1} )$ in this order, then $J_i \cup L_i$ is a closed curve for 
every $i=1\stb n$. Since the points $a$ and $b$ are outside the closed
disc $\cl U$, we have 
$$\sumin w(L_i ; a)=w(\del U ;a)=0 \quad \text{and} \quad 
\sumin w(L_i ; b)=w(\del U ;b)=0.$$
Thus
$$\sumin w(J_i \cup L_i ; a) =\sumin w(J_i ; a) +\sumin w(L_i ;a)=
w(J ;a) +0=\pm 2\pi$$
and
$$\sumin w(J_i \cup L_i ; b) =\sumin w(J_i ; b) +\sumin w(L_i ;b)=
w(J ;b) +0=0.$$
Therefore, we can choose an index $i$ such that $w(J_i \cup L_i ; a)
\ne w(J_i \cup L_i ; b).$

We prove that $\ga (t_{i-1} )$ and $\ga (t_i )$
belong to different connected components
of $D\se A$. Since $\ga (t_{i-1} ),$ $\ga (t_i ) \in \cl U$, this will prove 
the statement of the lemma.

Suppose that $\ga (t_{i-1} )$ and $\ga (t_i )$
belong to the same connected component of $D\se A$. Then 
there is an arc $L \su D\se A$ with endpoints 
$\ga (t_{i-1} )$ and $\ga (t_i )$. 

We have $w(L_i \cup L;a)= w(L_i \cup L;b)=0$, since
the closed curve $L_i \cup L$ is in $D$, while $a$ and $b$ are not in $D$.
Therefore, we have $w(L_i ;a)=- w(L ;a)=w(\ol L ;a)$ and 
$w(L_i ;b)=- w(L ;b)=w(\ol L ;b)$, where $\ol L$ is obtained from $L$ 
by changing its direction. Then we have
\begin{align*}
w(J_i \cup \ol L ;a)-w(J_i \cup \ol L ;b)& = w(J_i ;a)- w(L ;a)
-w(J_i ;b)+ w(L ;b) =\\
& = w(J_i ;a)+ w(L_i ;a)
-w(J_i ;b)- w(L_i ;b) =\\
&=w(J_i \cup L_i ;b) - w(J_i \cup L_i ;a)  \ne 0.
\end{align*}
Consequently, $a$ and $b$ belong to different connected components of the open set
$\sik \se (J_i \cup \ol L )$. This, however, is impossible, as $a,b\in A$,
$(J_i \cup \ol L )\cap A=\emp$, and $A$ is connected.
\hfill $\square$

\hfill \eject

\bsk 
\noi
Department of Mathematics,
University of Chicago 

\noi
Chicago, IL 60637, USA

\bsk 
\noi
Department of Analysis, Institute of Mathematics, E\"otv\"os Lor\'and University

\noi
Budapest, P\'azm\'any P\'eter s\'et\'any 1/C, 1117 Hungary

\msk \noi
{\rm e-mail:} csornyei@math.uchicago.edu (Marianna Cs\"ornyei),

\noi
herakornelia@gmail.com
(Korn\'elia H\'era), 

\noi
laczk@cs.elte.hu
(Mikl\'os Laczkovich)

\end{document}